# The BiEntropy of Some Knots on the Simple Cubic Lattice

v3.11

Grenville J. Croll

grenvillecroll@gmail.com

*Binary representations of the trefoil and other knots of up to ten crossings in the simple cubic lattice were created. The BiEntropy of each knot was computed using a variety of binary encodings and compared against controls. This showed that binary encoded knots are highly disordered information objects. The BiEntropy of knots on the simple cubic lattice increases slightly as the number of crossings and length of encoding increases. We show that the non-alternating knots of nine and ten crossings are more disordered than the alternating knots of nine and ten crossings.*

## 1 INTRODUCTION

We developed our BiEntropy function [Croll, 2013] as a simple method for determining the order and disorder of finite binary strings of arbitrary length. We successfully tested our BiEntropy function in a variety of domains including number theory, cryptography, quantitative finance and random number generation. We have started further work enumerating the algebras of BiEntropy [Croll, 2014] for possible application in bit string physics [Noyes, 1997]. Other teams have been able to successfully apply BiEntropy in cryptographic, internet information processing, mobile computing and random number generation domains [Costa et al, 2015][Jin & Zeng, 2015][Jin et al 2016][Kotě et al, 2014][Stakhanova et al, 2016].

We have become aware of Kauffman's work on knots [Kauffman, 2001] and physics [Kauffman, 2013] over a period of ten years or more. It was inevitable that we should at some point want to consider whether we could measure the BiEntropy of knots expressed in binary form.

Although the entropy of knots has been addressed from a variety of mathematical and statistical perspectives in the past [Franks & Williams, 1985][Baiesi, 2010] [van Rensburg & Rechnitzer, 2011] we believe our work is the first to quantify the entropy of a knot from its precise geometrical configuration using a generically applicable methodology.

This paper is based upon on the work of [Scharein et al, 2009] who have computed and tabulated the minimal or near minimal forms of knots of up to ten crossings on the simple cubic lattice. Using Scharein et al as a data source, we convert their NEWSUD knot encodings into various binary equivalents and compute the BiEntropy of the resultant binary encoded knots. In this paper we provide a brief introduction to Shannon's Entropy [Shannon, 1948], Binary Derivatives and define the several forms of BiEntropy we have previously used. Since knots are cyclical structures we describe a new variation of BiEntropy, Knot BiEntropy, which is effectively adapted to their geometric form. We describe the NEWSUD and subsequent binary encodings and then present our results in graphical and tabular form.

## 2 SHANNON ENTROPY, BINARY DERIVATIVES & WEIGHTING METHODS

### 2.1 Shannon Entropy

Shannon's Entropy of a binary string $s = s_1, \ldots, s_n$ where $P(s_i=1) = p$ (and $0 \log_2 0$ is defined to be 0) is:

$$H(p) = -p \log_2 p - (1-p) \log_2 (1-p)$$

For perfectly ordered strings which are all 1's or all 0's i.e. $p = 0$ or $p = 1$, $H(p)$ returns 0. Where $p = 0.5$, $H(p)$ returns 1, reflecting maximum variety. However, for a string such as 01010101, where $p = 0.5$, $H(p)$ also returns 1, ignoring completely the periodic nature of the string. BiEntropy seeks to compensate for this omission of consideration of the periodicity (i.e. the order and disorder) of a string by using the binary derivatives of *s*.

## 2.2 Binary Derivatives, Binary Knot Derivatives & Periodicity

In our previous work on the BiEntropy of linear strings, the first binary derivative of *s*, $d_1(s)$, is the binary string of length $n - 1$ formed by XORing adjacent pairs of digits. In the present work on the BiEntropy of knots, which are cyclical structures, the first binary knot derivative of *s*, $d_1(s)$, is the binary string of length *n* formed by XORing adjacent pairs of digits *including the last and the first*.

In either case we refer to the *k*th derivative of *s* $d_k(s)$ as the binary derivative of $d_{k-1}(s)$. There are $n-1$ binary derivatives of *s* in the linear string version of BiEntropy and infinitely many in the knot version of which we use the first *n-1*. In either case let $p(k)$ denote the observed fraction of 1's in $d_k(s)$ where $p(0)$ denotes the fraction of 1's in *s*.

By calculating the (first) *n-1* binary derivatives of *s* we can discover the existence of repetitive patterns in binary knots and strings of arbitrary (even) length. If a binary string is *periodic* then $d_{n-1}(s) = 0$. A binary string is *aperiodic* if $d_{n-1}(s) = 1$ or all 1's. A binary string is *nperiodic* if $d_{n-1}(s) = 0$, but it is not periodic.

For example, the first binary knot derivative of 01010101 (with *period*, $P = 2$) is 11111111 ($P = 1$), following which all the higher derivatives are all 0's indicating periodicity. The third derivative of 00010001 ($P = 4$) is 11111111, following which again all the higher derivatives are 0, again indicating periodicity. The seventh binary knot derivative of 10100011 is 11111111 indicating aperiodicity.

The properties of infinite binary strings, their derivatives and the notions of periodicity and eventual periodicity are given more fully in [Nathanson, 1971] and [Goka, 1970]. We rely here solely upon the binary derivatives of a *finite* string or *finite cyclical* knot to resolve the issue of the degree of periodicity (and hence the degree of order and disorder) within the string or knot.

## 2.3 Weighting Methods

In order to obtain a single value that represents the order and disorder of a binary string or knot we need to combine the Shannon Entropies of the derivatives in one or more ways. One obvious method is to combine them using the powers of two such that the last used binary derivative has the most weight (as the last derivative gives a final determination of periodicity and aperiodicity) and also in order that the influences of the weights of each derivative are arithmetically separated from each other. We abbreviate this method "BiEn" noting its meaning ("good") in the French language. Another method is to use a logarithmic weighting such that the order and disorder discovered in a binary string of *any* length through the enumeration of *all* its binary derivatives is fully taken into account, though in an exponentially diminishing fashion. We abbreviate this method Tres Bien or TBiEn for short. Linear weights would of course be LBiEn ("Les BiEn") and zero weights would be PBiEn ("Pas BiEn").

## 3 BIENTROPY

BiEntropy, or BiEn for short, is a weighted average of the Shannon binary entropies of the string and the first *n-2* binary derivatives of the string using a simple power law. This version of BiEntropy is suitable for shorter binary strings where $n \leq 32$ approximately as the weights of the first derivatives tend rapidly to zero.

$$\text{BiEn}(s) = \left(1 / (2^{n-1} - 1)\right) \left( \sum_{k=0}^{n-2} ((-p(k) \log_2 p(k) - (1 - p(k)) \log_2 (1 - p(k)))) 2^k \right)$$

The derivative $d_{n-1}$ is not used as it makes no contribution to the total entropy. The highest weight is assigned to the highest (i.e. last used) derivative $d_{n-2}$.

If the higher derivatives of an arbitrarily long binary string are periodic, then the whole sequence exhibits periodicity. For strings where the latter derivatives are not periodic, or for all strings in any case, we can use a second version of BiEntropy, which uses a Logarithmic weighting, to evaluate the complete set of a long series of binary derivatives.

$$\text{TBiEn}(s) = \left(1 / \sum_{k=0}^{n-2} \log_2 (k+2)\right) \left( \sum_{k=0}^{n-2} (-p(k) \log_2 p(k) - (1 - p(k)) \log_2 (1 - p(k))) \log_2 (k+2) \right)$$

The logarithmic weighting or (TBiEn for short) again gives greater weight to the higher derivatives

We illustrate in Table 1A the calculation in a spreadsheet of the logarithmic knot version of BiEntropy, which we abbreviate KTBiEn. We use KTBiEn exclusively in the analytics of this paper as the minimum length of a knot encoded in binary is 24 * 3 = 72 bits.

**Table 1A – Computing the Logarithmic Knot BiEntropy (KTBiEn) of an 8-bit string.**

| Binary Expansion of N | 1's | Len(n) | p | (1-p) | -p.log(p) | -(1-p).log(1-p) | BiEn | k | log(k+2) | BiEn.log(k+2) |
|---|---|---|---|---|---|---|---|---|---|---|
| 1 0 1 0 1 1 1 0 | 5 | 8 | 0.63 | 0.38 | 0.42 | 0.53 | 0.95 | 0 | 1.0 | 0.95 |
| 1 1 1 1 0 0 1 1 | 6 | 8 | 0.75 | 0.25 | 0.31 | 0.50 | 0.81 | 1 | 1.6 | 1.29 |
| 0 0 0 1 0 1 0 0 | 2 | 8 | 0.25 | 0.75 | 0.50 | 0.31 | 0.81 | 2 | 2.0 | 1.62 |
| 0 0 1 1 1 1 0 0 | 4 | 8 | 0.50 | 0.50 | 0.50 | 0.50 | 1.00 | 3 | 2.3 | 2.32 |
| 0 1 0 0 0 1 0 0 | 2 | 8 | 0.25 | 0.75 | 0.50 | 0.31 | 0.81 | 4 | 2.6 | 2.10 |
| 1 1 0 0 1 1 0 0 | 4 | 8 | 0.50 | 0.50 | 0.50 | 0.50 | 1.00 | 5 | 2.8 | 2.81 |
| 0 1 0 1 0 1 0 1 | 4 | 8 | 0.50 | 0.50 | 0.50 | 0.50 | 1.00 | 6 | 3.0 | 3.00 |
| 1 1 1 1 1 1 1 1 | | | | | | | | | | |
| | | | | | | | 6.39 | 21 | 15.2992 | 14.08924 |
| | | | | | | | | | Knot TBiEn ( s ) | **0.920913** |

We show in Table 1B a spreadsheet based calculation of the Logarithmic Knot BiEntropy (KTBiEn) of the 4 and 8 Bit strings. Note the symmetry about the diagonal and the consistently high values of BiEntropy for half of the strings.

**Table 1B – Logarithmic Knot BiEntropy (KTBiEn) of the 4 and 8-bit strings**

The BiEntropy of the perfectly ordered strings 00 and 11 is zero. The BiEntropy of the perfectly disordered strings 10 and 01 is one. For finite strings of arbitrary length greater than two, the BiEntropy is greater than or equal to zero and less than one, as there are no perfectly ordered strings of length greater than two.

**4 MINIMAL KNOTS IN THE SIMPLE CUBIC LATTICE**

[Diao, 1993] proved some time ago that the minimum number of steps needed to construct a knot on the simple cubic lattice was 24, and that the only knot available at this length is $3_1$, the trefoil knot - depicted in Figures 1 and 2. Minimum lengths for $4_1$ and $5_1$ of 30 and 34 steps respectively have been recently proven as outlined in Scharein et al. Except for these three proofs, determination of the minimum number of steps for knots on the cubic lattice is by heuristic methods and is therefore incomplete and somewhat uncertain.

**Figure 1 The Trefoil knot of canonical form 2-50 in the simple cubic lattice. The axis of symmetry is central and near vertical**

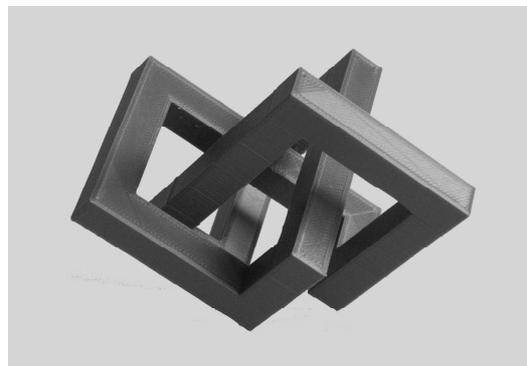

**Figure 2 Alternative View of Trefoil knot of canonical form 2-50. The axis of symmetry runs from bottom left to top right.**

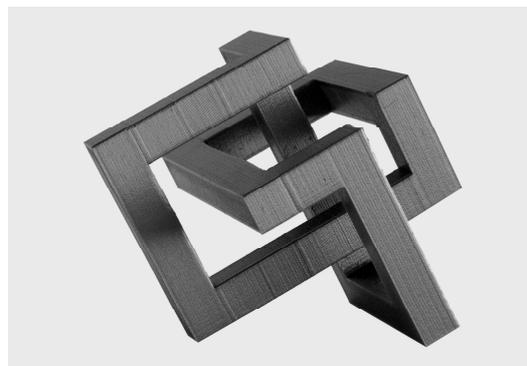

Scharein et al correct and complete earlier work to show by an exhaustive method that there are 75 minimum canonical forms of the trefoil knot (in three distinct groups) and tabulate each form in their Table 6 using the NEWSUD encoding. A NEWSUD string is a sequence of letters from the set {N, E, W, S, U, D} representing edge directions in the cubic lattice of North, East, West, South, Up, and Down. For example: DDDEEUUSWWWNNEEDSSSUUNNW is the trefoil 2-50 depicted in Figure 1 and 2 above. Scharein at al document the geometric realisation of each of the 75 canonical forms of the trefoil knot using stick and ball diagrams with the start point denoted by a small dot.

Scharein et al then give in their Table 7 a single example of a definite, potential or probable minimum encoding of each knot from $3_1$ up to $10_{165}$. Including the Unknot, the Granny and the Square knot by way of comparison. Other tables give the numbers of canonical forms of these more complex knots. Tables 6 & 7 of Scharein et al is the raw data for the rest of this paper.

## 5 BINARY ENCODINGS OF TREFOILS, KNOTS & CONTROLS

The simplest and most obvious encoding for the six member {N, E, W, S, U, D} set is a three bit encoding such as {000, 001, 010, 011, 100, 101}. Having allocated a binary encoding, it is trivial to compute the BiEntropy of the encoded knot. An encoding this simple would be naïve however as a number of issues need to be taken into consideration. We outline these issues in the following sections and introduce some of the thinking behind the design of the monte carlo experiments we perform.

### 5.1 Selection of encoding bits

There are 8! / (8-6)! = 20,160 ways of allocating a three bit binary number to each member of the NEWSUD set. Each encoding would result in a differing BiEntropy for each knot as the different sub patterns created in the resultant binary string representing the knot would be different. For example, the Unknot DEUW could be encoded in a three bit encoding as 000111010101 or 101010110011 which have differing BiEntropies. It is necessary therefore to compute the BiEntropy of a large enough sample of encodings in order to obtain statistically repeatable and reliable results.

### 5.2 Length of encoding

The existence and influence of sub patterns in the binary string representing a knot will vary according to the length of the encoding. For example, within an 8 bit encoding the Unknot DEUW could be represented as 00000000111111110000000011111110 which is relatively well ordered with a low BiEntropy or as 01110101110001101101001001001001 which is relatively disordered with a high BiEntropy. With a longer binary encoding, the three bits corresponding to the direction instructions implied in the NEWSUD encoding become lost in the noise. Note that there are 256! / 250! = 2.65E14 ways of selecting an 8 bit NEWSUD encoding. After some experimentation with 3 and 4 bit encodings we settled on an 8 bit encoding.

In our monte carlo analysis, we used two different encodings ENCODING_A and ENCODING_B. Each encoding consisted of 256 different randomly selected 8 bit encodings from the myriad available. For example, the first four encodings (out of 256) of ENCODING_A were:

**Table 2 – First four encodings (out of 256) of the NEWSUD characters of ENCODING_A**

| N | E | W | S | U | D |
|---|---|---|---|---|---|
| 44 | 82 | 201 | 21 | 245 | 214 |
| 32 | 231 | 121 | 87 | 252 | 179 |
| 6 | 28 | 80 | 33 | 137 | 77 |
| 111 | 52 | 163 | 120 | 238 | 147 |

**Table 3 – First four encodings (out of 256) of the NEWSUD characters of ENCODING_B**

| N | E | W | S | U | D |
|---|---|---|---|---|---|
| 84 | 41 | 102 | 101 | 67 | 222 |
| 43 | 107 | 20 | 118 | 66 | 113 |
| 227 | 111 | 65 | 189 | 142 | 99 |
| 145 | 108 | 159 | 3 | 248 | 240 |

Where the 44 of the N in the first row of ENCODING_A translates to 00101100 and the 3 of the S in the fourth row of ENCODING_B translates to 00000011. The unknot DEUW encoded using the first row of ENCODING_B would be {222; 41; 67; 102} in decimal or 11011110001010010100001101100110 in binary.

## 5.3 Start points

We changed the design of BiEntropy to accommodate the cyclical or hoop like nature of knots [Amson &Bowden, 2005]. The Knot version of BiEntropy detects the identity of the rotated versions of every eight bit string (e.g. 00000001, 00000010, …..10000000) and computes identical BiEntropies for each of the eight rotated versions. Thus we were able to use the NEWSUD strings exactly as given in Scharein et al for the trefoils and knots without needing to start the BiEntropy computation from a randomly allocated start point within each NEWSUD string. The elimination of a requirement for a start point following the decision to develop the knot version of BiEntropy lead to a duplication of most of the computations and results in this paper. The two sets of similar results provided a simple means of systematically comparing and checking our work across the two methods. We report only one set of results, for the Knot version of BiEntropy.

## 5.4 Trefoils and Knots

We worked with two data sets. Table 6 of Scharein et al documents the NEWSUD encodings for the 75 Canonical forms of the Trefoil knot. These results are denoted TREF_N. Table 7 of Scharein et al gives the NEWSUD encoding for a single suspected (or actual in three cases) minimal form of all the knots from $3_1$ through to $10_{165}$, which we denote SOME_N. Note that there was a mistake in $3_1$ of Table 7 as this was only 23 characters long. The initial "D" had been omitted due to a typographical error, which we corrected. We excluded the unknot, granny and square knots of Table 7 from our analysis giving 249 knots in total.

## 5.5 Randomised Trefoil and Knot controls

We computed the distribution of NEWSUD instructions within TREF_N and created two randomised controls TREF_A and TREF_B with similar NEWSUD character distributions. Each string in TREF_A and TREF_B is 24 characters long. TREF_A and TREF B contain 75 NEWSUD entries to match TREF_N.

We computed the distribution of NEWSUD instructions within SOME_N and created two randomised controls SOME_A and SOME_B with similar NEWSUD character distributions. Each string in SOME_A and SOME_B varied from 24 to 64 characters in length depending upon the knot length in the equivalent row of SOME_N. SOME_A and SOME_B contain 249 NEWSUD entries to match SOME_N.

The resulting four NEWSUD controls were pure random and were not computed to be self avoiding polygons. There is a very small probability that a control might contain a knot.

## 6 COMPUTATIONS

After a little experimentation to determine the sensitivities of various encoding lengths and so on we performed two similar computations in the 'C' language, firstly on the Trefoils and secondly on the knots of up to 10 crossings. These took ten minutes and one hour respectively on a 4GHz single core Pentium.

We examined every trefoil (in TREF_N) and every knot (in SOME_N) using each of the 256 encodings in firstly ENCODING_A and secondly ENCODING_B. For each trefoil, knot and encoding we created a binary string *s* and computed the logarithmic knot BiEntropy (KTBiEn(*s*)). We then repeated the exercise with the two pairs of control data TREF_A, TREF_B, SOME_A and SOME_B using both ENCODING_A and ENCODING_B.

Note that the 8 bit binary encoding of a trefoil knot encoded in NEWSUD is 24 * 8 = 192 bits long and that the 8 bit encoding of a more complex knot of up to 10 crossings is 64 * 8 = 512 bits in length.

Copyright © 2018 Grenville J. Croll. All Rights Reserved
6/11

# 7 RESULTS

The first result is that the BiEntropies of the Trefoils in TREF_N and the Knots in SOME_N computed using ENCODING_A and ENCODING_B, gave almost identical results across the two encodings. This was despite ENCODING_A and ENCODING_B being completely different from each other and having each been separately generated randomly from a very large space ($O_{E14}$). The strong correlations between the results for the two encodings are shown in Figures 3 and 4.

This result was almost exactly repeated in the random controls TREF_A, TREF_B, SOME_A and SOME_B using the two encodings. This suggests that we can reliably use sampled encodings in this and other experimental domains in order to obtain reliable estimates of BiEntropy. We have combined the results from ENCODING_A and ENCODING_B in the presentation of our remaining results.

**Figure 3 BiEntropy of Trefoil Knots (TREF_N) in the Simple Cubic Lattice**

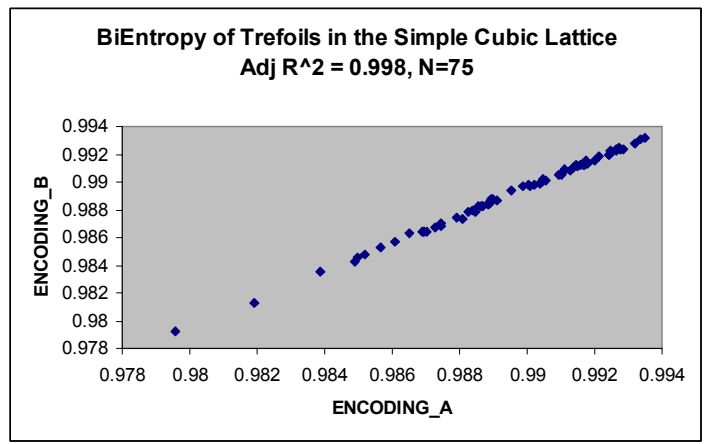

**Figure 4 BiEntropy of Knots up to 10 Crossings (SOME_N) in the Simple Cubic Lattice**

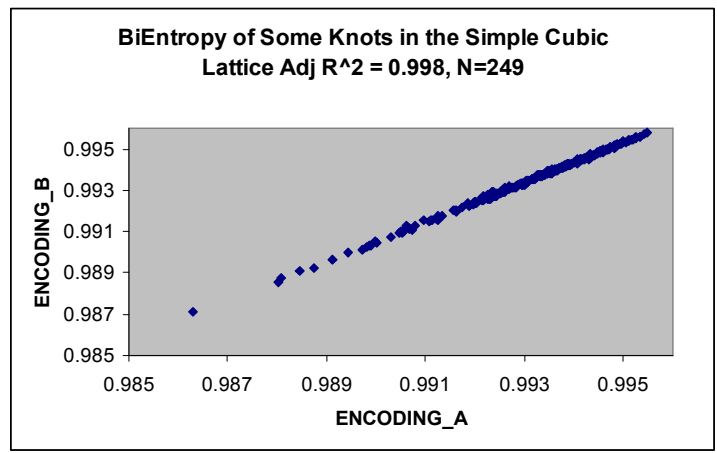

Secondly, as expected, the BiEntropy of knots on the simple cubic lattice was high (>0.979) indicating that these knots are relatively disordered binary objects. The knot with the lowest BiEntropy is the 2-50 variation of the trefoil which has a clear 3 way symmetry about a central axis, which can be determined from the photographs in Figures 1 & 2. We show the small variance in BiEntropy among the 75 canonical forms of the Trefoils in Figure 5 and the Knots in Figure 6.

Finally, the BiEntropies of the random controls TREF_A, TREF_B and SOME_A, SOME_B are slightly (0.3%) but significantly (p<0.001) higher than the corresponding BiEntropies of all of the trefoils and all of the knots. Thus indicating that the closed cyclical knots are more ordered forms.

**Figure 5 BiEntropy of the Trefoil Knots**

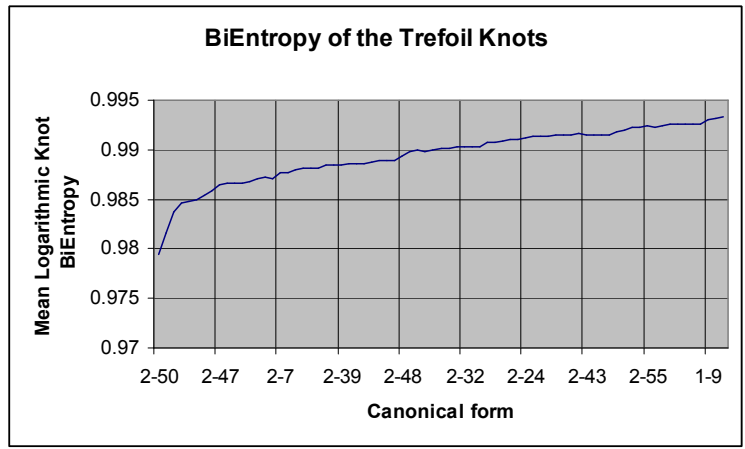

Figures 3 and 4 show in addition that whereas there is reasonable certainty as to which trefoils and knots have the lowest BiEntropy, there is sample size based uncertainty for the more complex knots. With larger sample sizes it may be possible to exactly sequence the knots in terms of their BiEntropy. The 8 bit encoding space is probably large enough to permit much more refined estimates of BiEntropy for the more complex knots, but at a significant computational cost.

**Figure 6 Mean BiEntropy of some knots up to 10 crossings**

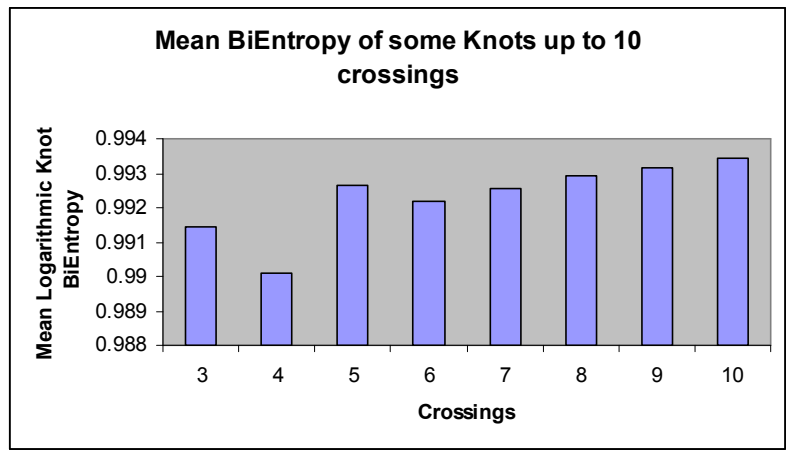

**Figure 7 Mean BiEntropy of some knots up to 10 crossings by NEWSUD encoding length**

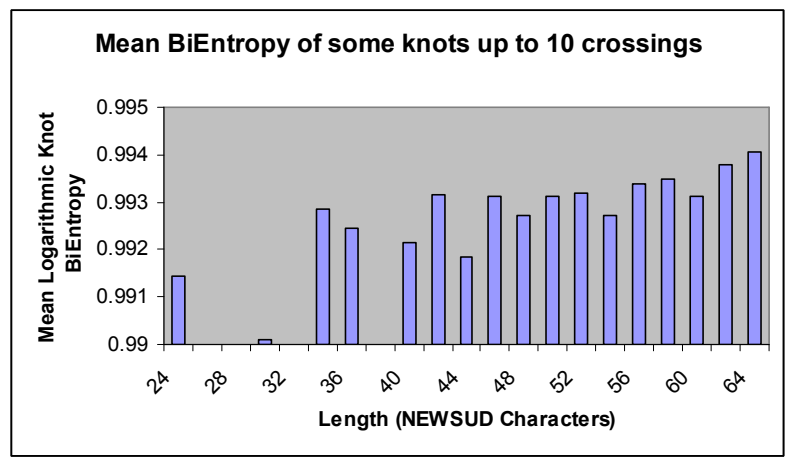

As expected, the mean BiEntropy of more complex knots with greater numbers of crossings slightly increased which we show in Figure 6, noting that there are single or few data points for the knots with lower crossings. As also expected, the BiEntropy of knots increased as the length of their NEWSUD encoding increased (Figure 7) noting again the low number of data points for the lower knots.

There were some small but significant (9 Crossings $p < 0.01$, 10 Crossings $p<0.001$) differences between the BiEntropies of the non-alternating (or quasi-alternating [Jablan, 2014]) versus alternating knots of differing lengths and crossings, which we show in Figures 8 and 9.

**Figure 8 The BiEntropy of Non-Alternating and Alternating knots of 9 crossings**

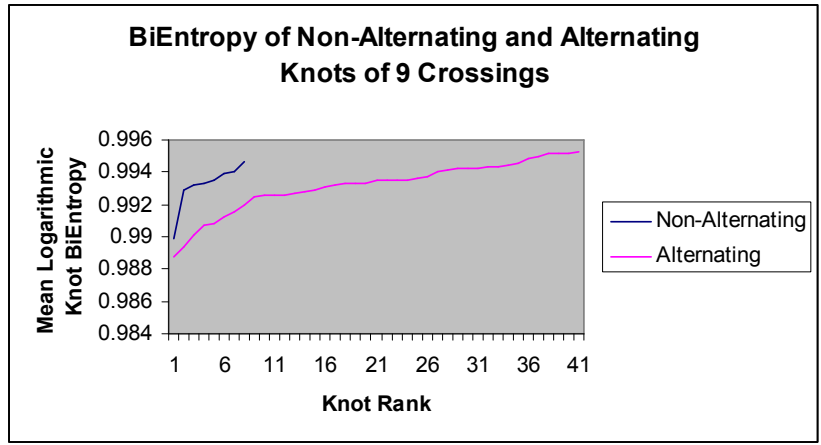

**Figure 9 The BiEntropy of Non-Alternating and Alternating knots of 10 crossings**

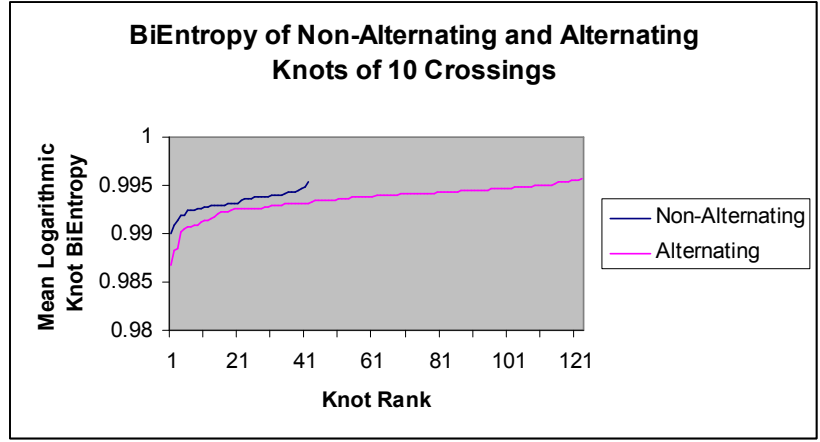

Note that the alternating and non-alternating knots are distinguished by their Dowker encodings. The integers of the Dowker encoding for an alternating knot are always positive. The integers of the Dowker encoding for a non-alternating knot are mixed positive and negative.

## 7 SUMMARY

We have adapted our BiEntropy measure for the cyclical world of knots. We have encoded some knots on the simple cubic lattice with an 8-bit binary coding and then evaluated the resulting long binary strings using BiEntropy.

More precisely, we have measured the Logarithmic Knot BiEntropy (KTBiEn) of the 75 canonical forms of the Trefoil knot and a single instance of each of the 249 distinct knots up to and including ten crossings.

Our results show that binary encoded knots are highly disordered binary objects. Despite the high level of disorder, we were able to discriminate between the entropy of knots of increasing crossings and length of encoding, between some alternating and non-alternating knots and between the trefoils of various canonical forms. We discovered that the 2-50 variation of the trefoil knot has the lowest BiEntropy of all of the knots we examined and exhibits a marked symmetry.

We have shown that we can use a simple monte carlo method to produce statistically consistent results in this domain with relatively low computational effort.

This study is constrained by our use of the limited output from a single prior experiment, but hopefully should point the way for wider studies into the entropy – the order and disorder – of binary encoded knots and other information structures.

## ACKNOWLEDGEMENTS

We thank Drs. John C. Amson (ANPA) and Simon Thorne (EuSpRIG) for their comments on an earlier draft of this paper and in particular Dr Amson's forthright insistence that we modify BiEntropy for the cyclical "bithoops" case. We thank Pascal Fleck for 3D knot design & printing services.## REFERENCES

Amson, J.C., & Bowden, K., (2005) A Combinatoric Bit-Hoop System. Proc. ANPA, Cambridge, "Against Bull", pp407-440

Baiesi, M., Orlandini, E., & Stella, A. L. (2010). The entropic cost to tie a knot. Journal of Statistical Mechanics: Theory and Experiment, 2010(06), P06012. http://arxiv.org/abs/1003.5134

Costa, R., Boccardo, D., Pirmez, L., & Rust, L. F. (2015) Hiding Cryptographic Keys of Embedded Systems. Advances in Information Science and Computer Engineering, ISBN: 978-1-61804-276-7

Croll, G. J. (2013). BiEntropy-The Approximate Entropy of a Finite Binary String. arXiv preprint arXiv:1305.0954.

Croll, G. J. (2014). BiEntropy–the Measurement and Algebras of Order and Disorder in Finite Binary Strings. Scientific Essays in Honor of H Pierre Noyes on the Occasion of His 90th Birthday. ISBN: 978-981-4579-36-0

Diao, Y. (1993). Minimal knotted polygons on the cubic lattice. Journal of Knot Theory and its Ramifications, 2(04), 413-425.

Franks, J., & Williams, R. F. (1985). Entropy and knots. Transactions of the American Mathematical Society, 241-253.

Goka, T., (1970) An operator on binary sequences, SIAM Rev., 12 (1970), pp. 264-266

Jablan, S. (2014). Tables of quasi-alternating knots with at most 12 crossings. arXiv preprint arXiv:1404.4965.

Jin, R., & Zeng, K. (2015, September). Physical layer key agreement under signal injection attacks. In Communications and Network Security (CNS), 2015 IEEE Conference on (pp. 254-262). IEEE.

Jin, R., Shi, L., Zeng, K., Pande, A., & Mohapatra, P. (2016). MagPairing: Pairing Smartphones in Close Proximity Using Magnetometers. IEEE Transactions on Information Forensics and Security, 11(6), 1306-1320.

Kauffman, L. H. (2001). Knots and physics (Vol. 1). World scientific. ISBN: 978-981-02-4111-7

Kauffman, L. H. (2013). Non-Commutative Worlds and Classical Constraints. Scientific Essays in Honor of H Pierre Noyes on the Occasion of His 90th Birthday, ISBN: 978-981-4579-36-0

Kotě, V., Molata, V., & Jakovenko, J. (2014). Enhanced Generic Architecture for Safety Increase of True Random Number Generators. Electroscope. 2014,

Nathanson, M. B. (1971). Derivatives of binary sequences. SIAM Journal on Applied Mathematics, 21(3), 407-412.

Noyes, H. P. (1997). A Short Introduction to Bit-String Physics. arXiv preprint hep-th/9707020.

van Rensburg, E. J., & Rechnitzer, A. (2011). Minimal knotted polygons in cubic lattices. Journal of Statistical Mechanics: Theory and Experiment, 2011(09), P09008.

Scharein, R., Ishihara, K., Arsuaga, J., Diao, Y., Shimokawa, K., & Vazquez, M. (2009). Bounds for the minimum step number of knots in the simple cubic lattice. Journal of Physics A: Mathematical and Theoretical, 42(47), 475006.

Shannon, C. E., (1948), A Mathematical Theory of Communication, Bell System Technical Journal 27 (3): 379–423

Stakhanova, N. (2016, March). An Entropy Based Encrypted Traffic Classifier. In Information and Communications Security: 17th International Conference, ICICS 2015, Beijing, China, December 9-11, 2015, Revised Selected Papers (Vol. 9543, p. 282). Springer.